\documentclass[12pt]{amsart}
\usepackage{amsmath, amssymb, amsthm, amscd}
\usepackage{graphics}

\newtheorem{theorem}      {Theorem}

\newtheorem{lemma}        {Lemma}

\newtheorem{proposition}  {Proposition}

\theoremstyle{definition}

\newtheorem*{notation}    {Notation}

\newcommand{\binomial}[2]  {{\left(\begin{array}{c}#1\\#2\end{array}\right)}}
\newcommand{\qalg}         {\mathbb Q^{\mathrm{a}}}

\newcommand{\norm}[1]	   {\ensuremath{{\, \left\| \, {#1}\, \right\|_2 \,}}}
\newcommand{\punct}[1]	   {\ \text{#1}}
\newcommand{\cond}[1]      {\ensuremath{\kappa_{\text{\tiny #1}}}}
\newcommand{\relgap}       {\nu}
\newcommand{\ratio}        {\rho}
\newcommand{\function}[5]  {\[ \begin{array}{rrcl}
#1:& #2 &\rightarrow & #3\\ & #4 & \mapsto & #5 \end{array}\]}
\title{Condition number bounds for problems with integer coefficients} 

\author{Gregorio Malajovich\footnote{
Departamento de Matem\'atica Aplicada, Universidade Federal
do Rio de Janeiro. Caixa Postal 68530, CEP 21945, Rio de Janeiro,
RJ, Brasil. e-mail: gregorio@labma.ufrj.br.
On leave at MSRI, 1000 Centennial Drive, Berkeley CA
94720-5070. e-mail: gregorio@msri.org
}}
\date{February 12, 1999}

\begin{document}
\begin{abstract}
An apriori bound for the condition number associated to 
each of the following problems is given: general linear 
equation solving, minimum squares, non-symmetric eigenvalue 
problems, solving univariate polynomials, solving systems 
of multivariate polynomials. It is assumed that the input 
has integer coefficients and is not on the degenerate locus 
of the respective problem (i.e. the condition number is 
finite). Then condition numbers are bounded in terms of the 
dimension and of the bit-size of the input.
\par
In the same setting, bounds are given for the speed of
convergence of the following iterative algorithms: QR
without shift for the symmetric eigenvalue problem, and
Graeffe iteration for univariate polynomials.
\end{abstract}
\maketitle

\section{Introduction}

In most of the numerical analysis literature, complexity
and stability of numerical algorithms are usually estimated
in terms of the problem instance dimension and of a `condition number'.
\par
For instance, the complexity of solving an $n \times n$ 
linear system $Ax=b$ is usually estimated in terms of the
dimension $n$ (actually the input size is $n(n+1)$) and of
the condition number $\cond{} (A) = \norm{A} \norm{A^{-1}}$.
\par
There is a set of problems instances with $\cond{}(A) = \infty$,
and in most cases it makes no sense to attempt solving those
problem instances. There are also problem instances (in our case, matrices)
close to the locus of degenerate problem instances. Those will have a
large condition number, and will be said to be {\em 
ill-conditioned}. 
\par
It is usually accepted that ill-conditioned
problem instances are hard to solve. Thus, for complexity purposes a
problem instance with a large condition number should be considered
`large'. Therefore, when considering problems defined for
real inputs, a reasonable measure for the input size
would be (in our example): $n^2 \log_2 \cond{}(A)$.
(Compare to ~\cite{SMALE97} Formula~2.1 and paragraph below. 
See also the discussion in ~\cite{BP}, Chapter~3, Section~1).
\medskip
\par
Another tradition, derived from classical complexity theory
and pervasive in several branches of literature (such as linear
programming), is to consider the subset of problems instances with integer
coefficients. Hence the input size is the number of coefficients
times the bit-size of the largest coefficient (in absolute value).
\par
In this paper, the following classical problems of numerical
analysis are considered:
\newcounter{numitems}
\begin{enumerate}
\item Solving a general $n \times n$ system of linear equations.
\item Minimal squares problem for a full-rank matrix.
\item Non-symmetric eigenvalue problem.
\item Solution of one univariate polynomial.
\item Solution of a non-degenerate system of $n$ polynomial equations
     in $n$ variables.
\setcounter{numitems}{\value{enumi}}
\end{enumerate}
\par
All those problems share the feature mentioned above: there is a 
degenerate locus, and problem instances with real coefficients can be as close
to the degenerate locus as wished. This implies that they can be
arbitrarily ill-conditioned.
\par
However, in Theorems~\ref{th1} to ~\ref{th5} below, we provide
bounds for the condition number of problems instances with
integer coefficients and not in the degenerate locus. 
Those bounds depend on the dimension (size) of the problem instance and on the
bit-size of its coefficients. 
\medskip
\par
In the analysis of iterative algorithms, one further considers a
certain quantity that can be used to bound the speed of convergence
and hence the number of iterations
to obtain a given approximation. For instance, for
power methods (or QR iteration without shift)
in the symmetric eigenvalue problem, one can bound 
the number of steps in terms of the desired accuracy and of the
ratio between different eigenvalues. The farther this number is
from 1, the faster is the convergence.
\par
Once again, if input has real coefficients, this quantity can be
arbitrarily close to 1. However, explicit bounds for that quantity
will be given for inputs with integer coefficients for
\begin{enumerate}
\setcounter{enumi}{\value{numitems}}
\item QR iteration without shift for the Symmetric Eigenvalue Problem.
\item Graeffe iteration for solving univariate polynomials.
\end{enumerate}
\medskip
\par
The reader should be warned that the results herein are worst
case estimates, and are overly pessimistic for application purposes.
The main motivation for those results is to convert numerical
analysis estimates into `polynomial time' estimates, not the
opposite.

\section{Statement of main results}

\begin{notation} \norm{.} stands for the {\em 2-norm}: if
$x \in \mathbb R^n$ or $\mathbb C^n$, then 
\[ \norm{x} = \sqrt{\sum_{i=1}^n |x_i|^2} 
\punct{.}
\]
If $A$ is a matrix, then
\[
\norm{A} = \sup_{\norm{x} = 1} \norm{Ax}
\punct{.}
\]
\end{notation}

\subsection{Linear equation solving}

The first problem considered is linear equation solving:
given an $n \times n$ matrix $A$ and a vector $b \in
\mathbb R^n$, find $x \in \mathbb R^n$ such that
$ Ax = b $.
\par
Its condition number (with respect to the 2-norm) is
defined as
\[
\cond{}(A) = \norm{A} \norm{A^{-1}}
\punct{.}
\]
\par
Comprehensive treatment of the perturbation theory for this
problem can be found in the literature, such as~\cite{DEMMEL}
Section~2.2, ~\cite{HIGHAM} Chapter~7, ~\cite{TB} Lecture~12,
etc...

\begin{theorem}\label{th1}
Let $A$ be an $n \times n$ matrix with integer
coefficients. If $A$ is invertible, then
\[
\cond{}(A) \le n^{\frac{n}{2}+1}
\left(
\max_{i,j} |A_{ij}|
\right) ^n
\punct{.}
\]
\end{theorem}

\par
No originality is claimed for Theorem~\ref{th1}. This result is
included for completeness and because its proof is elementary,
yet illustrates the principle behind the other results.

\subsection{Minimal squares}

The second problem in the list is minimal squares fitting. 
Let $A$ be an $m \times n$ matrix, $m \ge n$, with full rank,
and let $b \in \mathbb R^m$. One has to find $x$ to minimize
$\norm{A x - b}^2$. 
\par
Let $r = Ax - b$ be the residual, we are minimizing $\norm{r}^2$.
Let 
\[
\sin \theta = \frac{\norm{r}}{\norm{b}}
\punct{.}
\]
\par
According to~\cite{DEMMEL} p. 117  
(Compare to ~\cite{TB} Lecture~18 and~\cite{HIGHAM}
Section~19.1),
the condition number of the linear least squares problem
is
\[
\cond{LS} (A,b)
=
\frac{2 \cond{}(A)}{\cos \theta} + \tan \theta \ \cond{}(A)^2
\punct{.}
\]
\par
Since we do not assume $A$ to be square, we need to give a
new definition for $\cond{}(A)$. Let $\sigma_{\text{\tiny MAX}}(A)$
and $\sigma_{\text{\tiny MIN}}(A)$ be respectively the largest and
the smallest singular values of $A$. Then set
\[
\cond{}(A) = \frac{\sigma_{\text{\tiny MAX}}(A)}{\sigma_{\text{\tiny MIN}}(A)}
\punct{.}
\]
When $m=n$, this definition is equal to the previous one.
\medskip
\par
The singular locus is now the set of pairs $(A,b)$ such that
$A$ does not have full rank (i.e. $\sigma_{\text{\tiny MIN}}(A)=0$)
or such that $\norm{r} = \norm{b}$ (i.e. $b$ is orthogonal
to the image of $A$).
\par
The result is:
\begin{theorem}\label{th2}
Let $A$ be an $m \times n$ matrix with integer coefficients,
and assume that $A$ has full rank. Let $b \in \mathbb Z^m$.
Set $H = \max _{i,j} \left( |A_{ij}|, |b_i| \right)$. Then
if $b$ is not orthogonal to the image of $A$, we have:
\[
\cond{LS} (A,b) \le
3n ^{\frac{n}{2}+1} m^{n+\frac{1}{2}}H^{2n+1}
\punct{.}
\]
\end{theorem}
\par

\subsection{Non-symmetric eigenvalue problem}

Let $A$ be an $n \times n$ matrix and let $\lambda$
be a single eigenvalue of $A$. The condition number
of $\lambda$ depends on the angle
between the left and right eigenvectors:
\par
Let $x$, $y$ be respectively right and left norm-1
eigenvectors of $A$ associated to $\lambda$:
$Ax = \lambda x$, $y^* A = \lambda y^*$, and
$\norm{x} = \norm{y} = 1$. Then
\[
\cond{NSE}(A, \lambda) = \sec(\widehat{x,y}) = \frac{1}{y^* x}
\punct{.}
\]
\par
See ~\cite{DEMMEL} Theorem~4.4 p.~149 for references.
\begin{theorem}\label{th3}
Let $A$ be an $n \times n$ matrix with integer
coefficients, and let $\lambda$ be a single
eigenvalue of $A$. Then
\[
\cond{NSE}(A, \lambda) 
\le
n^{3n} 2^{2n} \left( 2 \sqrt{n} H(A) \right)^{2 n^3 - 2 n}
\punct{.}
\]
\end{theorem}

\subsection{Solving univariate polynomials}

The condition number (in affine space) for solving a univariate
polynomial $f(x) = \sum_{i=0}^d f_i x^i$ can
be defined~(\cite{BCSS} page 228) as:
\[
\mu(f) = \max_{\zeta \in \mathbb C: f(\zeta)=0} \mu(f,\zeta)
\punct{,}
\]
where
\[
\mu(f, \zeta) = 
\frac { \left( \sum_{i=0}^d |\zeta|^{2i} \right)^{\frac{1}{2}}} 
      { |f'(\zeta)|}
\punct{.}
\]
\par
The degenerate locus is the set of polynomials with a multiple
root or with a root at infinity.

\begin{theorem}\label{th4}
Let $f: x \mapsto \sum_{i=0}^d f_i x^i$ be a univariate polynomial with
integer coefficients, without multiple roots. Then
\[
\mu(f) \le 2^{2d^2 - 2} d^{2d} \left(\max |f_i| \right)^{2d^2}
\punct{.}
\]
\end{theorem}

\subsection{Solving systems of polynomials}

A similar condition number exists for systems of polynomials.
However, for the purpose of condition number theory, it is usually 
convenient to homogenize the equations and to study the
perturbation theory of the `roots' in complex projective
space. This can also be seen as a change of metric, that
simplifies the formula of the condition number and of several
theorems (See~\cite{BCSS} Chapters~10, 12, 13).
\par
Let $f = (f_1,$ $\cdots,$ $f_n)$ be a system of polynomials in
variables $x_1, \cdots, x_n$. We homogenize the system by
multiplying each coefficient 
$f_{iJ} x^J = f_{iJ} x_1^{J_1} x_2^{J_2} \cdots x_n^{J_n}$ of $f_i$ by 
$x_0^{J_0}$, where we choose $J_0 = \deg f_i - (J_1 + \cdots J_n)$.
We obtain a system of homogeneous polynomials in $n+1$ variables,
that we call $F=(F_1, \cdots, F_n)$. The natural space for the
roots of $F$ is projective space $\mathbb P^n$, defined as the
space of all `rays' 
\[
(x_0: \cdots : x_n) =
\{ (\lambda x_0, \cdots, \lambda x_n): \lambda \in \mathbb C\}
\punct{.}
\]
\par
where $x_0$, \dots, $x_n$ are not all equal to $0$.
\par
Every finite root $(x_1, \cdots, x_n)$ of $f$ corresponds to the
projective root of $F$ given by $(1: x_1: \cdots: x_n)$, and 
projective roots of $F$ correspond either to a finite root
of $f$ or to a root `at infinity'.
\par
Suppose that the coefficients of $f$ (hence of $F$) are made
to depend upon a parameter $t$. The condition number 
bounds the absolute speed of the roots of $F$ (in
projective space) with respect to the absolute speed 
of the coefficients of $F$. Recall that the roots $\zeta$
of $F$ are in projective space, so their speed vector $\dot \zeta$
belongs to the tangent space $T_{\zeta} \mathbb P^n$. 
\par
The condition number of $F$ at a root 
turns out to be:
\[
\mu(F,\zeta) = \norm{F} \norm{ \left( \mathrm{D}F(\zeta)_{| T_\zeta}\right) ^{-1}
\left[ \begin{matrix}
  \norm{\zeta}^{d_1 -1}\\
  & \ddots\\
  & & \norm{\zeta}^{d_n -1}\\
\end{matrix} \right] }
\]
where $\zeta \in \mathbb C^{n+1}$ is such that
$(\zeta_0 : \cdots : \zeta_n)$ is a root of $F$ 
(See Proposition 7c in Page 230 of~\cite{BCSS}).
We did not define the norm of a polynomial yet. Above, 
$\norm{.}$ stands for the unitary invariant
norm~(See \cite{WEYL} Chapter III-7 or ~\cite{BCSS} Section~12.1), 
that is the most reasonable
generalization of the 2-norm to spaces of
polynomials: 
\begin{notation}
Let $G$ be a homogeneous degree $d$ polynomial in $n+1$ variables.
Then 
\[
  \norm{G} = \sqrt{ \sum_J \frac{ |G_J|^2 }{\binomial{d}{J}} } 
\]
where $\binomial{d}{J}$ is $\frac{d!}{J_0 ! \cdots J_n!}$.
Let $F$ be a system of homogeneous polynomials. Then
\[
  \norm{F} = \sqrt{ \sum \norm{F_i}^2 }
\punct{.}
\]
\end{notation}
\par
With these definitions, the number $\mu(F,\zeta)$ is invariant
under scalings of $F$, $\zeta$, and under the action of the
unitary group $U(n+1)$, where an element $Q \in SU(n+1)$ acts
by $Q:(F,\zeta) \mapsto (F \circ Q, Q \zeta)$.
\par 
In order to define the condition number of a system of
$n$ equations in $n$ variables, we set:
\[
\mu(f) = \max_{\zeta} \mu(F,\zeta)
\]
where $\zeta$ ranges over the roots of $F$. (Another 
possibility is to restrict $\zeta$ to the non-degenerate
roots of $F$. This would make no difference in this paper).
The following theorem is true if one restricts $\zeta$ to
any subset of the roots of $F$.

\begin{theorem}\label{th5}
Let $f$ be a system of $n$ polynomial equations in $n$ variables,
with integer coefficients. We write $H(f)$ for the maximum of the absolute
value of the coefficients of $f$, $S(f)$ for the number of
non-zero coefficients of $f$ and $D$ for $\max d_i$.
Assume that $\mu(f)$ is finite. Then
\[
\mu(f) \le ((n+1) S H(f))^{D^{cn}}
\]
where $c$ is an universal constant.

\end{theorem}

\subsection{Symmetric eigenvalue problem}
\label{sep}
\par
Let $A$ be an $n \times n$ real positive symmetric matrix, and let
$\lambda_1 \ge \lambda_2 \ge \cdots \lambda_n \ge 0$ be its
eigenvalues. 
\par
Unlike the non-symmetric eigenvalue problem,
the symmetric eigenvalue problem has absolute condition number
always equal to 1 (See \cite{DEMMEL} Theorem~5.1. See also ~cite{PARLETT} 
Fact~1.11 p.16).
\par
However, when using an iterative algorithm, the
ratio of eigenvalues
\[
\ratio(A) = \min_{j > i} \frac{\lambda_j}{\lambda_i} 
\]
may play an important role for estimating convergence. For instance,
according to ~\cite{TB} Theorem~28.4, the QR algorithm without shift
converges linearly with speed $\frac{1}{\ratio(A)}$. Convergence may
get slower when $\ratio(A) \rightarrow 1$. Therefore one can bound
the speed of convergence by bounding
\[
\relgap(A) = 
\ratio(A) - 1 
=
\min_{j > i} \frac{\lambda_j - \lambda_i}{\lambda_i} 
\]
above from zero. If $\relgap(A) > \delta_0$, then
$\ratio(A) > 1+\delta_0$.
After $k > \lceil \frac{1}{\delta_0} \rceil$ iterations,
one gets
\[
\ratio(A)^{k} > 1 + k \delta_0 \ge 2 
\punct{.}
\]
\par
Thus it suffices to perform 
$O( \frac{1}{\delta_0} \log_2 \frac{1}{\delta_1})$
iterations to obtain a result with accuracy $\delta_1$.
\par
Also, the quantity 
$\relgap(A)^{-1}$ can also be interpreted as a condition number
for the eigenvectors (See ~\cite{DEMMEL} Theorem~5.7 p.~208).
We will show here that
\par

\begin{theorem}\label{th6}
Let $A$ be an $n \times n$ matrix with integer coefficients.
Then
\[
\relgap(A)^{-1} \ge 8^{-n} (4n)^{-n^2} \left( \max_{i,j} |A_{ij}| \right)^{-2n^2} 
\punct{.}
\]
\end{theorem}

\subsection{Graeffe iteration}

Let $f: x \mapsto \sum_{i=0}^d f_i x^i$, $f_d=1$
be a monic univariate polynomial
with zeros $\zeta_1, \cdots, \zeta_d$. Those zeros can be ordered
such that
\[
|\zeta_1 | \ge |\zeta_2| \ge \cdots \ge |\zeta_d|
\punct{.}
\]
\par
The Graeffe operator maps the polynomial $f(x) = \prod_{i=1}^d (x-\zeta_i)$
into the polynomial $Gf(x) = (-1)^d f(\sqrt{x}) f(-\sqrt{x}) =
\prod_{i=1}^d (x-\zeta_i^2)$. 
\par
In ~\cite{OSTROWSKI,OSTROWSKI2}, it is explained how to
recover the actual roots of $f$ after a certain number of 
Graeffe iterations, with a good approximation. The number of required
iterations depends on the ratio:
\[
\ratio(f) = \max_{|\zeta_j| > |\zeta_i|} \frac{|\zeta_j|}{|\zeta_i|}
\punct{.}
\]
\par
Unlike in Section~\ref{sep}, we do not require here
that the roots have different absolute value. 
We consider also the auxiliary quantity
\[
\relgap(f) = \ratio(f) - 1 =  \max_{|\zeta_j| > |\zeta_i|}
\frac{|\zeta_j| - |\zeta_i|}{|\zeta_i|}
\punct{.}
\]
\par
By the above definitions, the `condition number' 
$\relgap(f)^{-1}$ is always finite. In order to
recover the roots within relative precision 
$\delta$, the number
of Graeffe iterations to perform is
\[
O(\log \relgap(f)^{-1} + \log d + \log \log \delta^{-1})
\punct{.}
\]
\par
For clarity of exposition, we will show that
bound under a special hypothesis: all the roots should be different 
positive real numbers. For the general case, see ~\cite{MZ97} and
~\cite{MZ98}. Also, all estimates here are `up to the first order',
and quadratic error terms will be discarded.
\medskip
\par
After $k$ steps of Graeffe iteration one
obtains the polynomial
\[
g(x) = G^kf(x) = 
\sum_{i=0}^d g_i x^i = \prod_{i=1}^d (x - \zeta_i^{2^k})
\]
with $\ratio(g) = \ratio(f)^{2^k}$. 
\par
Expanding each $g_i$ as the $(d-i)$-th elementary symmetric function
of the $\zeta_i^{2^k}$, one obtains
\begin{eqnarray*}
g_0     &=& \sigma_d ( \zeta_1^{2^k}, \cdots ,  \zeta_d^{2^k} )\\
g_1     &=& \sigma_{d-1} ( \zeta_1^{2^k}, \cdots ,  \zeta_d^{2^k} )\\
&\vdots&\\
g_{d-1} &=& \sigma_1 ( \zeta_1^{2^k}, \cdots ,  \zeta_d^{2^k} )\\
g_d     &=& 1
\end{eqnarray*}
\par
We can use the special hypothesis to bound
\[
\left( \zeta_1 \zeta_2 \cdots \zeta_{d-i} \right) ^{2^k} = g_i(1+\delta_i)
\]
with $|\delta_i| < \frac{2^d}{\ratio(g)} + \text{h.o.t.}$ 
Hence
\[
\zeta_i^{2^k} = \frac{g_{d-i+1}}{g_{d-i}} (1 + \delta_i')
\]
with $|\delta_i'| < \frac{2^{d+1}}{\ratio(g)} + \text{h.o.t.}$
\par
Since we assumed the $\zeta_i$ are all positive, we can recover
them by taking $2^k$-th roots
\[
\zeta_i = \left(\frac{g_{d-i+1}}{g_{d-i}}\right)^{2^{-k}} (1 + \delta_i'')
\punct{.}
\]
\par
with $|\delta_i''| \le \frac{2^{d+1-k}}{\ratio(g)} + \text{h.o.t.}$
\par
Now we can use the estimate on $\ratio(g) = \ratio(f)^{2^k}$ to
deduce that
$O(\log \relgap(f)^{-1} + \log \delta^{-1})$ steps are sufficient
to obtain a relative
precision $\delta$ in the roots. Indeed after $k_1 =  \log_2 \relgap(f)^{-1}$
steps,
\[
\ratio(G^{k_1}f) = \ratio(f)^{2 ^{ \log_2 \relgap(f)^{-1}}} \ge 2 
\punct{.}
\]
\par
After extra $k_2 = \log_2 (d +1+ \log_2 \delta^{-1})$ steps, one gets
\[
\ratio(G^{k_1+k_2}(f)) > 2^{2^{\log_2 (d+1+\log_2 \delta^{-1})}} 
= 2^d \delta^{-1} 
\punct{.}
\]
\par
So we can set $k = k_1 + k_2 + 1$, the last $1$ to get rid of the
high order terms, to deduce that $|\delta_i''| < \delta$.
\medskip
\par
\begin{theorem}\label{th7}
Let $f: x \mapsto \sum_{i=0}^d f_i x^i$ be a polynomial with
integer coefficients. Then
\[
\relgap(f)^{-1} > \left( 8 \max |f_i| \right)^{-2d}
\punct{.}
\]
\end{theorem}

\medskip
\par
This says that Graeffe iteration is `polynomial time', in the sense
that we can obtain relative accuracy $\delta$ of the roots after
\[
O( d \log \max |f_i| + \log \log \delta^{-1})
\]
steps.

\section{Background material}

The proof of Theorems~\ref{th3} to ~\ref{th7} will make
use of the {\em absolute multiplicative height function} $H$
to bound inequalities involving
algebraic numbers. 
\par
The construction of the height function $H$ is quite standard in number
theory and we refer the reader to ~\cite{LANG} Chapter~II 
or to~\cite{SILVERMAN}
pages 205--214. For applications to complexity theory, see~\cite{BCSS}
Chapter~7 and ~\cite{M98}. 
\par
The height function is naturally defined in the projectivization
$\mathbb P^n(\qalg)$ of the algebraic numbers $\qalg$. It returns a 
real number $\ge 1$. We can also extend it to complex projective space 
$\mathbb P^n$ by setting $H(P) = \infty$ when $P \not \in \mathbb P^n(\qalg)$.
We will adopt this convention in order to simplify the notation of
domains and ranges.
\par
Also, if $x = (x_1, \cdots, x_n) \in \mathbb C^n$, we can define its
height as $H(x) = H(x_1: \cdots : x_n:1)$. 
\par
We can also define the height of matrices, polynomials and systems of
polynomials as the height of the vector of all the coefficients.
\medskip
\par
The following properties of heights will be used in the sequel.
First of all, we can explicitly write the height of a vector 
with integer coefficients as:

\begin{proposition}\label{integers}
If $u \in \mathbb Z^n$, then $H(u) = \max_{1\le i \le n} |u_i|$,
where $|\, . \,|$ is the standard absolute value.
\end{proposition}

Proposition~\ref{integers} follows from the construction of 
the height function. One immediate consequence is that if
$v \in \mathbb Q^n$, then $H(v) = \max |m v_i|, |m|$ where
$m$ is the greatest common denominator of the $v_i$'s. 
\par
We can use the following fact to bound the height of the
roots of an integral polynomial:

\begin{proposition}\label{roots}
If $f: x \mapsto f(x) = \sum_{i=0}^d f_i x^i$ is a non-zero
polynomial with integer coefficients, and if $x$ is a root of $f$,
then $H(x) \le 2 \max |f_i|$.
\end{proposition}
\par
Proposition~\ref{roots} is Theorem~5 in ~\cite{M98}. Compare
with Theorem~5.9 in~\cite{SILVERMAN}, where the coefficients of 
$f$ are algebraic numbers.
\par
We can use a bound on the height to bound absolute values above and below:

\begin{proposition}\label{bound}
Let $K$ be an algebraic extension of $\mathbb Q$,
and let $x \in K$, $x \ne 0$. Then
\[
H(x)^{-\deg[ K:\mathbb Q]} \le |x| \le H(x)^{\deg[ K:\mathbb Q]}
\punct{.}
\]
\end{proposition}

The height of a vector and of its coordinates can be related by:

\begin{proposition}\label{coordinate}
\[
H(x_1) \le H(x_1, \cdots, x_n) \le H(x_1) H(x_2) \cdots H(x_n)
\punct{.}
\]
\end{proposition}

Propositions~\ref{bound} and~\ref{coordinate} follow immediately
from the construction of the height function.
The height function is invariant under permutation of coordinates,
and also:

\begin{proposition}\label{conjugate}
	Let $K$ be an algebraic extension of $\mathbb Q$, and
	let $g \in \text{Gal}_{[K:\mathbb Q]}$. Then for any $x \in K$,
	$H(g(x)) = H(x)$
\end{proposition}

\par
Proposition~\ref{conjugate} is Lemma~5.10 in~\cite{SILVERMAN}.

\begin{proposition}\label{morphisms}
Let 
\function{F=(F_0, \cdots, F_m)}
{\mathbb C^{n_1} \times \mathbb C^{n_2} \times \cdots \times \mathbb C^{n_k}}
{\mathbb C^{m+1}}
{P^1, \cdots , P^k}
{F(P^1, \cdots, P^k)}
be a system of multi-homogeneous polynomials with algebraic coefficients, 
where each $F_i$ has degree $d_j$ in variables $P^j$.  
Let the $P^j$ be algebraic.  Then 
\[ H(F(P)) \le (\max S(f_i)) H(F) H(P^1)^{d_1} \cdots H(P^k)^{d_k} 
\punct{.}
\] 
\end{proposition}

\par
In the case $k=1$, this is similar to Theorem~5.6 in~\cite{SILVERMAN} 
(where $\max S(f_i)$ is not given explicitly). For the general case 
see Theorem~4 in~\cite{M98}.

\begin{proposition}\label{polynomials}
Let
\function
{G=(G_1, \cdots, G_m)}
{\mathbb C^{n_1} \times \mathbb C^{n_2} \times \cdots 
\times \mathbb C^{n_k}}
{\mathbb C^{m}}
{Q^1, \cdots , Q^k}
{G(Q^1, \cdots, Q^k)}
be a system of polynomials with algebraic
coefficients, where each $G_i$ has degree at most $d_j$ in
variables $Q^j$. Let the $Q^j$ be algebraic. Then
\[
H(G(Q)) \le (\max S(G_i)) H(G)
H(Q^1)^{d_1} \cdots H(Q^k)^{d_k}
\punct{.}
\]
\end{proposition}

This is Corollary~1 in~\cite{M98}. Some consequences of this
are that $H(\sum_{i=1}^n x_i) \le n \prod H(x_i)$ and that
$H(\prod_{i=1}^n x_i) \le \prod H(x_i)$.

The following fact follows also from the construction
of heights:

\begin{proposition}\label{square}
If $x$ is an algebraic number,
\[
H(x^2) = H(x)^2
\punct{.}
\]
\end{proposition}

\par
Also, it makes sense to bound the height of the roots of
a system of polynomials with respect to the height, size
and degree of the system. Corollary~6 in ~\cite{KP} is:
\begin{proposition}\label{krickpardo}
[Krick and Pardo]
Let $f_1, \cdots , f_r$ $r \le n$ be polynomials
in $\mathbb Z[x_1, \cdots$ , $x_r]$ of degree and height 
bounded by $d \ge n$ and $\eta$, respectively,
and let $V$ denote the algebraic affine variety defined
by: $V = \{ x: f_1(x) = \cdots f_r(x) = 0 \}$.
\par
Then $V$ has at most $d^n$ isolated points, and their height verifies:
\[
\log_2 H(P) \in d^{O(n)} (\log_2 r + \log_2 \eta) 
\punct{.}
\]
\end{proposition}

\section{Proof of Theorems}

\begin{notation}
If $A$ is a real (resp. complex) matrix, then $A^*$ is the
real (resp. complex) transpose of $A$, $(A^*)_{ij}= \bar{A_{ji}}$.
The same convention will be used for vectors.
\par
The vectors of the canonical basis will be denoted by $e_1 =
[1,0,0,\cdots]^*$, $e_2 = [0,1,0,\cdots]^*$, etc... 
\end{notation}

\subsection{Proof of Theorem~\ref{th1}}

\begin{eqnarray*}
\norm{A} &=& \sup_{\norm{u}=1} \norm{A u} \text{\ by definition} \\
         &\le& \sum_j |u_j| \norm{ [A_{1j}, \cdots, A_{nj}]^* }
         \text{\ by triangular inequality}\\
         &\le& \sqrt{n} \max_j \norm{ [A_{1j}, \cdots, A_{nj}]^* }
	 \text{\ since \norm{u}=1}\\
	 &\le& n \max_{ij} |A_{ij}|
\end{eqnarray*}

\par
Let $A(i, u)$ be the matrix obtained by replacing 
the $i$-th column of $A$ by the vector $u$. Then if
$v = A^{-1} u$, Cramer's rule is:
\[
v_i = \frac{ \det A(i, u) } {\det A}
\punct{.}
\]
\par
Since $A$ has integer coefficients and
$\det A \ne 0$, one can always bound $|v_i| \le |\det A(i, u)|$.
By Hadamard inequality, this implies:
\begin{eqnarray*}
|v_i| &\le& \norm{u} \max_j \norm{ [A_{1j}, \cdots, A_{nj}]^* }^{n-1}\\
    &\le& \norm{u} (\sqrt{n})^{n-1} \left(\max_{i,j} |A_{ij}|\right)^{n-1}
\end{eqnarray*}
Therefore,
\begin{eqnarray*}
\norm{A^{-1}} &=& \sup_{\norm{u}=1} \norm{A^{-1} u} \text{\ by definition} \\
&\le& n^{\frac{n}{2}}  \left(\max_{i,j} |A_{ij}|\right)^{n-1} \\
\end{eqnarray*}
\par
Combining the bounds for \norm{A} and \norm{A^{-1}}, one obtains:
\[
\cond{}(A) \le n^{\frac{n}{2}+1}  \left(\max_{i,j} |A_{ij}|\right)^{n}
\punct{.}
\]

\subsection{Proof of Theorem~\ref{th2}}

In order to estimate $\cond{}(A)$, we write
\begin{eqnarray*}
\cond{}(A) &=& \sqrt{ \cond{} (A^* A) } \\
&\le & n^{\frac{n}{4}+\frac{1}{2}} H^{n} m^{n/2}
\end{eqnarray*}
\par
In order to bound $\cos \theta$, we use the assumption that
$b$ is not orthogonal to the image of $A$. Hence $\norm{A^* b} \ge 1$
and the `normal equation' $ A^* A x = A^* b $ implies:
\[
\norm{A^* A x} \ge 1
\punct{.}
\]
Therefore,
\[
\cos \theta = \frac{\norm{A^* A x}}{\norm{b}} \ge \frac{1}{\norm{b}}
\ge \frac{1}{H \sqrt{m}}
\]
and $\frac{1}{\cos \theta}$ and $\tan \theta$ are bounded above by
$H \sqrt{m}$. Putting all together, 
\begin{eqnarray*}
\cond{LS} (A,b) &\le&
2 n ^{\frac{n}{4}+\frac{1}{2}} m^{\frac{n}{2}+\frac{1}{2}} H^{n+1}
+
n ^{\frac{n}{2}+1} m^{n+\frac{1}{2}} H^{2n+1}
\\
&\le&
3 n ^{\frac{n}{2}+1} m^{n+\frac{1}{2}} H^{2n+1}
\end{eqnarray*}

\subsection{Proof of Theorem~\ref{th3}}

\begin{lemma}\label{lemma1}
Let $B$ be an $n \times n$ matrix with integer
coefficients. Let $p(t) = \det (B - tI) = \sum p_i t^i$.
Then 
\[
\max |p_i| \le \left( 2 \sqrt{n} \max_{i,j} |B_{ij}| \right)^n
\punct{.}
\]
\end{lemma}

\begin{proof}[Proof of Lemma~\ref{lemma1}]
\[
p_i = \sum_C \pm \det C
\]
where $C$ ranges over the $(n-i)\times(n-i)$ sub-matrices of $B$ of
the form $C_{kl}=B{s_ks_l}$ for some $1 \le s_1 < \cdots < s_{n-1} \le n$. 
Hence
\begin{eqnarray*}
|p_i| &\le& \binomial{n}{i} \max_C |\det C| \\
      &\le& \binomial{n}{i} \left( \sqrt{n-i} \max_{ij} |C_{ij}| \right)^{n-i}\\
      &\le& \left(2 \sqrt{n} \max_{i,j} |B_{ij}| \right)^n
\end{eqnarray*}
\end{proof}

\begin{lemma}\label{lemma2}
Let $A$ be an $n \times n$ matrix with integer coefficients and
let $\lambda$ be an eigenvalue of $A$. Then 
\[
H(\lambda) \le 2 \left(2 \sqrt{n} \max_{i,j}|A_{ij}| \right)^n
\punct{.}
\]
\end{lemma}

\begin{proof}[Proof of Lemma~\ref{lemma2}]
Apply Proposition~\ref{roots} to the polynomial $p(t)$
from Lemma~\ref{lemma1}.
\end{proof}

\begin{lemma}\label{lemma3}
Let $B$ be an $n \times n$ matrix with integer coefficients.
Let $q(t) = \det (B-tI + t e_n^* e_n) = \sum q_i t^i$. Then
\[
\max |q_i| \le 2 \left( 2 \sqrt{n} \max_{i,j} |B_{ij}| \right)^n
\punct{.}
\]
\end{lemma}

\begin{proof}[Proof of Lemma~\ref{lemma3}]
Let $p(t) = \det (B-tI)$ and let $r(t) = \det (\tilde B - t I)$
where $\tilde B$ is the $(n-1)\times(n-1)$ matrix obtained by deleting
the $n$-th row and the $n$-th column of $B$. Then, by multi-linearity
of the determinant,
\[
p(t) = q(t) \pm t\ r(t)
\punct{,}
\]
hence
\[
q(t) = p(t) \pm t\ r(t)
\punct{.}
\]
\par
Therefore, 
\begin{eqnarray*}
\max |q_i| &\le& \max |p_i| + \max |r_i| \\
&\le& 
\left(2 \sqrt{n} \max_{i,j} |B_{ij}| \right)^n
+
\left(2 \sqrt{n-1} \max_{i,j} |B_{ij}| \right)^{n-1}
\text{\ (Lemma~\ref{lemma1})}
\\
&\le&
2 \left(2 \sqrt{n} \max_{i,j} |B_{ij}| \right)^n
\end{eqnarray*}
\end{proof}

\begin{lemma}\label{lemma4}
Let $A$ be an $n\times n$ matrix with integer coefficients.
Let $\lambda$ be an isolated eigenvalue of $A$ and let
$x$ be an eigenvector associated to $\lambda$, $Ax = \lambda x$.
Then 
\[
H(x_1 : \cdots : x_n) \le 
n 2^n \left( 2 \sqrt {n} \max_{i,j}|A_{ij}| \right)^{n^2 -1}
\punct{.}
\]
\end{lemma}

\begin{proof}[Proof of Lemma~\ref{lemma4}]
Assume without loss of generality that the first $n-1$ lines
of $A-\lambda I$ are independent. Let $M_1$, \dots, $M_i$, \dots, $M_n$
be the sub-matrices obtained from $A-\lambda I$ by deleting the
last line and the $i$-th column. Then we can scale $x$ in such a
way that
\[
x_i = \pm \det M_i 
\punct{.}
\]
\par
We have $M_n = B_n - \lambda I$. By reordering rows and columns, we
obtain for each $i<n$ that $M_i$ is of the form:
\[
M_i = B_i - \lambda I + \lambda e_{n-1}^* e_{n-1}
\]
where $B_i$ is the sub-matrix of $A$ obtained by deleting the last
line and the $i$-th column.
\par
Set $q^{(i)} (\lambda) = \det M_i = \sum q^{(i)}_j \lambda^j$. Now by
Lemma~\ref{lemma3},
\[
\max |q^{(i)}_j| \le 2 \left( 2 \sqrt{n-1} \max_{k,l} |A_{kl}| \right)^{n-1}
\punct{.}
\]
\medskip
\par
We consider now the morphism:
\function{q}{\mathbb P}{\mathbb P^n}
{(\lambda:1)} { (q^{(1)}(\lambda): \cdots : q^{(n)}(\lambda)) }
\par
Then $x=q(\lambda)$ and
\begin{eqnarray*}
H(x) &=& H(q(\lambda)) \\
     &\le& n H(q) H(\lambda)^{n-1} 
     \\
     &\le& 
     n 2 \left( 2 \sqrt{n-1} \max_{i,j} |A_{ij}|\right)^{n-1}
     2^{n-1} \left( 2 \sqrt{n} \max_{i,j} |A_{ij}|\right)^{n(n-1)}
     \\
     &\le& 
     n 2^{n} \left( 2 \sqrt{n} \max_{i,j} |A_{ij}|\right)^{n^2 -1}
\end{eqnarray*}
	the first inequality because of  Proposition~\ref{morphisms},
	and the second because of  Lemma~\ref{lemma2}.
\end{proof}

\medskip
\par
\begin{proof}[End of the Proof of Theorem~\ref{th3}]
Proposition~\ref{polynomials} implies
\[
H(y^* x) \le n H(x) H(y)
\punct{.}
\]
\par
We claim that $\deg [ \mathbb Q[y^* x] : \mathbb Q] \le n$. Indeed,
$x$ and $y$ can be obtained by solving systems of linear equations
with coefficients in $\mathbb Q[\lambda]$, thus 
$x_i, y_i \in \mathbb Q[\lambda]$. Also, $\overline{y_i} 
\in \mathbb Q[\bar{\lambda}] = \mathbb Q[\lambda]$ so
$\deg [ \mathbb Q[y^* x] : \mathbb Q] \le n$ as claimed.
\par
By hypothesis $y^* x \ne 0$. Hence, by Proposition~\ref{bound},
\begin{eqnarray*}
|y^* x| &\ge& \left(n H(x) H(y) \right)^{-n} \\
&\ge& n^{-3n} 2^{-2n} \left( 2 \sqrt{n} \max_{i,j} |A_{ij}| \right)^{-2n^3+2n}
\end{eqnarray*}
\end{proof}
\subsection{Proof of Theorem~\ref{th4}}
\begin{proof}
  According to Proposition~\ref{roots}, 
\[
H(\zeta) \le 2 \max |f_i|
\punct{.}
\]
\par
Also,
\[
H(f'_i) \le d H(f) = d \max |f_i|
\]
and according to Proposition~\ref{polynomials}
\begin{eqnarray*}
H(f'(\zeta)) &\le& d H(f') H(\zeta)^{d-1} \\
             &\le& d^2 2^{d-1} H(f)^d  \\
\end{eqnarray*}
\par
and hence $|f'(\zeta)| \ge d^{-2d} H(f)^{-d^2}
2^{-d(d-1)}$. On the other hand,
\begin{eqnarray*}
\left(
\sum_{i=0}^d |\zeta|^{2i} 
\right)^\frac{1}{2}
&\le&
\sqrt{d+1}
\max(1,|\zeta|^d) 
\\
&\le&
\sqrt{d+1}
H(\zeta)^{d^2} 
\\
&\le&
\sqrt{d+1}
2^{d^2} H(f)^{d^2} 
\end{eqnarray*}
\par
Hence 
\[
\mu(f,\zeta) \le 2^{2d^2 -2} H(f)^{2 d^2} d^{2d}
\punct{.}
\]
\end{proof}

\subsection{Proof of Theorem~\ref{th5}}

\begin{lemma}\label{lemma5}
Let $A$ be an $n \times n$ invertible matrix with algebraic
coefficients. Then
\[
	H(A^{-1}) \le n H(A)^n
\punct{.}
\]
\end{lemma}

\begin{proof}[Proof of Lemma~\ref{lemma5}]
Let $A(i,j)$ be the sub-matrix of $A$ obtained by deleting
the $i$-th row and the $j$-th column. By Cramer's rule,
$\left(A^{-1}\right)_{ji} = \frac{\det A(i,j)}{\det A}$. Therefore we
should define the degree $n$ morphism:
\function{\varphi}{\mathbb P^{n^2}}{\mathbb P^{n^2}}
{\begin{array}{l}\left(A_{11}: A_{12}: \cdots \right.\\
\hspace{1cm} \left. \cdots : A_{nn} : 1\right)\end{array}}
{\begin{array}{l}\left( \det A(1,1): \det A(1,2): \cdots \right. \\
\hspace{1cm} \left. \cdots : \det A(n,n) : \det A\right)\end{array}}
\par
Then by Proposition~\ref{morphisms}, 
\begin{eqnarray*}
H( \varphi(A) ) &\le& n! H(A)^n H(\varphi) \\
&\le& n! H(A)^n
\end{eqnarray*}
\end{proof}

Let us fix the notations
\[
M
=
 \left[ \begin{matrix}
 \norm{\zeta}^{1-d_1}\\
 & \ddots\\
 & & \norm{\zeta}^{1-d_n}\\
 & & & \norm{\zeta}^{-1}\\
 \end{matrix} \right]
 \left[
 \begin{matrix}
 \mathrm{D}F(\zeta) \\
 \zeta^*
 \end{matrix}
 \right]
\]
and
\[
C =
 \mathrm{D}F(\zeta)_{| T_\zeta} ^{-1}
  \left[ \begin{matrix}
   \norm{\zeta}^{d_1-1}\\
    & \ddots\\
     & & \norm{\zeta}^{d_n-1}\\
      \end{matrix} \right]
\punct{.}
\]
\par
Let $\zeta \in \mathbb C^{n+1}$ be a fixed
representative for a root of $F$. 
Any $u \in T_{\zeta} \mathbb P^n$ can be written
as a vector in $\mathbb C^{n+1}$, orthogonal to
$\zeta$. Computing $u = Cv$ is the same as
solving $Mu = \left[ \begin{matrix}v\\0 \end{matrix}\right]$.
The operator $C$ 
is the same as $(M^{-1})_{|x_{n+1}=0}$. Therefore,
\[
\norm{C} \le \norm{M^{-1}}
\punct{.}
\]

\begin{lemma} \label{lemma6}
	In the conditions of Theorem~\ref{th5},
\[
H \left( M
 \right) 
 \le S \sqrt{n+1}^{D-1} H(\zeta)^{2 D - 2} D H(f) 
\punct{.}
\]
\end{lemma}

\begin{proof}[Proof of Lemma~\ref{lemma6}]
We apply Proposition~\ref{polynomials} to the
system:
\[
\zeta, N \mapsto 
\left( \cdots, N^{d_i -1} \frac{\partial F_i}{\partial x_j} (\zeta), 
\cdots, \bar \zeta_j, \cdots \right)
\]
with $N = \norm{\zeta}^{-1}$ to obtain
\[
H(M) \le
S H(\zeta)^{D-1} H(N)^{D-1} H(\mathrm{D}F)
\punct{.}
\]
\par
We can bound $H(\mathrm{D}F) \le D H(F)$ and $H(N) = H( \norm{\zeta})
= \sqrt{H(\sum |\zeta_i|^2)}$. We can apply Proposition~\ref{morphisms}
to the map
\function{\varphi}{\mathbb C^{n+1}}{\mathbb C}{\zeta, \bar \zeta}
{\sum \zeta_i, \bar \zeta_i}
to get $H(N^2) \le (n+1) H(\zeta) H(\bar \zeta)$. Proposition~\ref{conjugate}
implies $H(\zeta) = H(\bar \zeta)$, hence:
\[
H(N^2) \le (n+1) H(\zeta)^2
\]
and by Proposition~\ref{square},
\[
H(N) \le \sqrt{n+1} H(\zeta)
\punct{.}
\]
Thus, we can estimate that
\[
H(M) \le S  \sqrt{n+1}^{D-1} H(\zeta)^{2D-2} D H(f)
\punct{.}
\]
\end{proof}

\begin{proof}[End of the Proof of Theorem~\ref{th5}]
By definition of the norm, $\norm{f} \le S H(f)$. By 
Lemma~\ref{lemma5} and Lemma~\ref{lemma6}, 
we have:
\begin{eqnarray*}
H(C)
 &\le&
 (n+1) H(M) ^n \\
 &\le&
 (n+1) S^n \sqrt{n+1}^{nD-n} H(\zeta)^{2 nD - 2n} D^n H(f)^n 
 \\
\end{eqnarray*}
\par
Knowing that $\deg \left[\mathbb Q[\zeta] : \mathbb Q \right] \le D^n$,
we can use Proposition~\ref{bound} to deduce that
\begin{eqnarray*}
\norm{C} &\le& (n+1) H(C) \\
         &\le& (n+1)
\left( 
 (n+1) S^n \sqrt{n+1}^{nD-n} H(\zeta)^{2 nD - 2n} D^n H(f)^n 
\right)^{D^n}
\end{eqnarray*}
According to Proposition~\ref{krickpardo},
\[
H(\zeta) \le n H(f)^{D^{c'n}}
\]
\par
where $c'$ is a universal constant. Thus,
\begin{eqnarray*}
\mu(f,\zeta)
&\le&
S H(f) 
 (n+1)\left( (n+1) S^n \sqrt{n+1}^{nD-n} 
 n^{2 nD -2n} \right. \\
& & \hspace{4cm} \left. H(f)^{D^{c'n} (2 nD - 2n)} D^n H(f)^n\right)^{D^n} \\
&\le& 
((n+1) S H(f))^{D^{cn}}
\end{eqnarray*}
\par
where $c$ is a universal constant.

\end{proof}

\subsection{Proof of Theorem~\ref{th6}}

Lemma~\ref{lemma2} implies: 
\[
H(\lambda) \le 2 \left( 2 \sqrt{n} \max_{i,j} |A_{ij} \right) ^n
\punct{.}
\]
Hence (Proposition~\ref{polynomials}),
\[
H\left( \frac{\lambda_j}{\lambda_i} - 1\right) 
\le
8 (4n)^n \left( \max_{i,j} |A_{ij}| \right) ^{2n}
\punct{.}
\]
\par
Thus, by Proposition~\ref{bound},
\[
\left| \frac{\lambda_j}{\lambda_i} - 1\right| 
\le
8^n (4n)^{n^2} \left( \max_{i,j} |A_{ij}| \right) ^{2n^2}
\punct{.}
\]

\subsection{Proof of Theorem~\ref{th7}}

According to Proposition~\ref{roots}, 
\begin{eqnarray*}
H(\zeta_i) &\le& 2 H(f) \\
H(\zeta_j) &\le& 2 H(f) 
\end{eqnarray*}
\par
Moreover, $H(|\zeta_i|) \le H(\zeta_i)$ because $|\zeta_i|^2 = \zeta_i 
\bar{\zeta_i}$ and $H(\zeta_i) = H(\bar{\zeta_i})$ 
(Propositions~\ref{conjugate},~\ref{polynomials} and~\ref{square}). Thus,
\begin{eqnarray*}
H \left( \frac{ |\zeta_j| } {|\zeta_i|} -1 \right)
&\le& 2 H(|\zeta_i|) H(|\zeta_j|) \\
&\le& 8 H(f)^2
\end{eqnarray*}
\par
It follows that
\[
\relgap(f)^{-1} \ge \left( 8 H(f) \right) ^{-2\deg [\mathbb Q[|\zeta_i|,
|\zeta_j|]: \mathbb Q]} \ge \left( 8 H(f) \right) ^{-2d}
\punct{.}
\]

\section{Further comments}

As mentioned before, a reasonable definition for the
`real complexity' input size is the number  
of coefficients of a given problem instance, times the logarithm
of its condition number.
\par
Theorems~\ref{th1} to~\ref{th4} show that the `real complexity' input
size is no worse than a polynomial of the `classical complexity' input
size, for problem instances with integer coefficients. Theorem~\ref{th5}
also, if one considers $D^n$ as part of the input size. It may be 
possible to replace $D^n$ by the B\'ezout number $\prod d_i$,
that is the number of solutions of a generic system
of polynomials. 
\par
Since the `real complexity' of the problems considered can be bound by
common numerical analysis techniques, those Theorems provide a scheme
to convert `real complexity' bounds into `classical complexity' bounds.
\medskip
\par
The same idea is behind Theorems  ~\ref{th6} and ~\ref{th7}. In the case
of the iterative algorithms considered, the number of iterations for 
obtaining a certain approximation can also be bounded in terms of 
a `condition number'. In the case of problem instances with integer
coefficients, the `condition number' is also polynomially bounded in
terms of the input size.
\medskip
\par
Those Theorems have many features in common, and this is not a 
coincidence. A more general approach is to interpret the condition number 
as the inverse of the distance to the degenerate locus. This can
be bounded in terms of the height of the problem instance, and in
terms of the degenerate locus (degree, dimension, height). However,
bound obtained this way will be no sharper and possibly worse than
the direct bounds obtained by using the exact expression for the
condition number.
\medskip
\par
This paper was written while the author was visiting MSRI at Berkeley.
He wishes to thank MSRI for its generous support. 
Thanks to Bernard Deconinck, Jennifer Roveno, Paul Gross, Raquel,
and very special thanks to Paulo Ney de Souza and family.
\bibliographystyle{plain}
\bibliography{worstcond}

\end{document}